\documentclass[a4paper,12pt]{article}
\usepackage{comment}
\usepackage{cite}
\usepackage{amsmath,amssymb,amsfonts}
\usepackage[T1]{fontenc}
\usepackage[utf8]{inputenc}
\usepackage{graphicx}
\usepackage{fancyhdr}
\usepackage{float}
\usepackage{xcolor}
\usepackage{authblk}
\usepackage{mathrsfs}
\usepackage{empheq}
\usepackage[hyphens]{url}
\usepackage{hyperref}
\usepackage{amsthm}
\usepackage{tikz}
\usetikzlibrary{decorations.markings}

\usepackage{geometry}
\theoremstyle{plain}

\theoremstyle{definition}

\theoremstyle{remark}

\begin{document}

\title{Integral and arithmetic structures of alternating (zigzag) numbers $A_n$}

\author[$\dagger$]{Jean-Christophe {\sc Pain}\footnote{jean-christophe.pain@cea.fr}\\
\small
$^1$CEA, DAM, DIF, F-91297 Arpajon, France\\
$^2$Université Paris-Saclay, CEA, Laboratoire Matière en Conditions Extrêmes,\\ 
F-91680 Bruyères-le-Châtel, France
}

\date{}

\maketitle

\begin{abstract}
The alternating (zigzag) numbers $A_n$, counting the ascending alternating permutations of $\left\{1,\cdots,n\right\}$ and defined by the exponential generating function $\tan x+\sec x$, admit several classical combinatorial and analytic representations. In this work we unify and extend three complementary structures of $A_n$. First, starting from the Stirling number expansion of zigzag numbers, we derive a contour integral representation, as well as a positive Laplace–type integral representation
\[
A_n = 2^n \int_0^\infty e^{-y} f_n(y)\, dy, \qquad 
f_n(y) := \sum_{k=0}^{n} (-1)^k S(n,k) \left(\frac{y}{2}\right)^k,
\]
where the kernel \(f_n(y)\) is the polynomial generating function of Stirling numbers. A continuous interpolation of the discrete product (falling factorial) is introduced subsequently. This provides a direct analytic bridge between set partitions and Laplace asymptotics. Second, using the partial fraction expansion of $\tan$, we obtain the well-known hyperbolic integral representation
\[
A_{2n+1}=\frac{1}{\pi}\int_0^\infty\frac{y^{2n+1}}{\sinh(y/2)}\,dy,
\]
equivalently expressed in classical $\cosh$ form for $A_{2n}$. This representation interprets zigzag numbers as spectral moments associated with half–integer poles. The connection with Fourier analysis and Mellin transforms is also outlined. Finally, combining spectral expansions with Stirling identities, we derive congruence relations modulo primes for $A_n$. These results exhibit a dual analytic–combinatorial structure of zigzag numbers, linking partition expansions, trigonometric spectra, and arithmetic properties.
\end{abstract}

\section{Introduction}

Alternating permutations, also called zigzag permutations, have a rich historical background in combinatorics. They were first systematically studied by Désiré André in the late 19th century, who investigated the enumeration of permutations \(\sigma \in \mathfrak{S}_n\) satisfying the up-down pattern \(\sigma_1<\sigma_2>\sigma_3<\sigma_4>\cdots\) and established connections with the secant and tangent functions~\cite{Andre1881}. These permutations later appeared in the work of Euler and other mathematicians studying the interplay between combinatorial structures and trigonometric series, leading to the classical alternating numbers \(A_n\) as coefficients in the expansion of \(\tan x + \sec x\). Over time, alternating permutations have been linked to a variety of topics, including set partitions, Euler numbers, Bernoulli numbers, and more recently, analytic and arithmetic structures.

These combinatorial structures have been systematically studied in the classical works of Comtet~\cite{Comtet1974} and Stanley~\cite{Stanley2011,StanleyEC2}, and they appear naturally in connection with analytic combinatorics~\cite{Flajolet2009}. Early investigations of finite difference calculus and related sequences can be found in Nörlund~\cite{Norlund1924}.

The alternating (zigzag) numbers $A_n$, counting permutations $\sigma\in \mathfrak{S}_n$ with pattern
\[
\sigma_1<\sigma_2>\sigma_3<\sigma_4>\cdots,
\]
i.e., ascending alternating permutations, are characterized by the exponential generating function
\[
\tan x+\sec x=\sum_{n\ge0}A_n\frac{x^n}{n!}.
\]
They admit several classical representations involving Euler numbers, Stirling numbers, and trigonometric expansions. The zigzag numbers split naturally into even and odd subsequences. The even terms coincide with the absolute Euler (secant) numbers,
\[
A_{2n}=|E_{2n}|,
\]
where $E_{2n}$ are defined by
\[
\sec x=\sum_{n\ge0}(-1)^n E_{2n}\frac{x^{2n}}{(2n)!}.
\]
The odd terms are related instead to Bernoulli numbers,
\[
A_{2n+1}=(-1)^n\frac{2^{2n+2}(2^{2n+2}-1)}{2n+2}\,B_{2n+2},
\]
via the classical expansion of $\tan x$.

The purpose of this article is to unify and extend three complementary viewpoints. First, we present a derivation of a positive Laplace-type integral for $A_n$ starting from their Stirling number expansion. Second, we show that the hyperbolic integral form for $A_n$ arises naturally from partial fractions of the tangent function. Third, we investigate arithmetic consequences of these expansions and provide new congruence relations modulo primes. The first 10 alternating (zigzag) numbers $A_n$ are given in table~\ref{tab:zigzag_numbers}.

\begin{table}[H]
\centering
\caption{First 10 alternating (zigzag) numbers $A_n$}
\begin{tabular}{c|cccccccccc}
$n$ & 1 & 2 & 3 & 4 & 5 & 6 & 7 & 8 & 9 & 10 \\
\hline
$A_n$ & 1 & 1 & 2 & 5 & 16 & 61 & 272 & 1385 & 7936 & 50521
\end{tabular}
\label{tab:zigzag_numbers}
\end{table}

Section~2 presents an integral representation of the zigzag numbers $A_n$ derived from Stirling numbers, obtained by expressing the combinatorial sum over set partitions as a contour integral. This approach highlights the combinatorial structure of $A_n$ and provides a natural starting point for analytic methods, including asymptotic expansions. In Section~3, the same Stirling-number expansion is used to derive a strict Laplace-type integral, replacing factorials with their continuous integral representations. This yields a positive integral over $[0,\infty)$, well-suited for asymptotic analysis and continuous interpolation, and emphasizes the direct link between set partitions and Laplace integrals. Section~4 develops a hyperbolic integral representation of $A_n$ from the partial fraction decomposition of the tangent function, offering a spectral perspective complementary to the Stirling–Laplace representation. This form interprets the zigzag numbers as moments associated with half-integer poles and connects the trigonometric generating functions to classical analytic transforms. Section~5 explores these analytic connections, showing how the hyperbolic integral representations relate to Fourier coefficients of $\sec x + \tan x$ and to classical sums over half-integers via the Mellin transform. Finally, Section~6 studies arithmetic properties of $A_n$ based on the Stirling-number expansions, establishing congruences modulo primes and highlighting the role of set partitions in the arithmetic behavior of the zigzag numbers.

\section{Contour integral representation of zigzag numbers via Stirling numbers}

The zigzag (up/down) numbers $A_n$ admit exact expansions in terms of Stirling numbers of the second kind $S(n,k)$ (also denoted as $\left\{\!\begin{matrix} n \\ k \end{matrix}\!\right\}$). A classical identity (see Ref.~\cite{Comtet1974}) is
\begin{equation}
A_n = \sum_{k=0}^{n} (-1)^k\, 2^{\,n-k}\, k!\, S(n,k).
\label{eq:zigzag-stirling}
\end{equation}
We now combine this identity with the standard contour–integral representation of Stirling numbers (see Ref.~\cite{Flajolet2009}):
\begin{equation}
S(n,k)=\left\{\!\begin{matrix} n \\ k \end{matrix}\!\right\}
=
\frac{1}{2\pi i}
\oint_{\mathcal C}
\frac{(e^z-1)^k}{k!\, z^{n+1}}\,dz,
\label{eq:stirling-integral}
\end{equation}
where $\mathcal C$ is any small positively oriented contour around the origin. For large $n$, the term 
$$
\left(-\frac{e^z-1}{2}\right)^{n+1}
$$ 
is exponentially small near the saddle point $z_0 \sim 0$, because 
$$|e^z-1|/2 \ll 1$$ 
in this neighborhood. Therefore, its contribution to the integral is negligible compared to the leading term, allowing us to approximate
\[
A_n \sim \frac{2^{\,n+1}}{2\pi i} \oint_{\mathcal C} \frac{dz}{(e^z+1)\, z^{\,n+1}}.
\]
Substituting Eq.~\eqref{eq:stirling-integral} into Eq.~\eqref{eq:zigzag-stirling} and cancelling $k!$ gives
\begin{equation*}
A_n
=
\frac{1}{2\pi i}
\oint_{\mathcal C}
\frac{1}{z^{n+1}}
\sum_{k=0}^{n} (-1)^k 2^{\,n-k}(e^z-1)^k
\,dz.
\end{equation*}
The finite geometric sum evaluates to
\[
\sum_{k=0}^{n} (-1)^k 2^{\,n-k}(e^z-1)^k
=
2^n \sum_{k=0}^{n} 
\left(-\frac{e^z-1}{2}\right)^k,
\]
and hence we have
\begin{equation*}
A_n
=
\frac{2^n}{2\pi i}
\oint_{\mathcal C}
\frac{\displaystyle
\sum_{k=0}^{n} 
\left(-\tfrac{e^z-1}{2}\right)^k
}{z^{n+1}}
\,dz.
\end{equation*}
Using the closed form of a finite geometric series,
\[
\sum_{k=0}^{n} r^k=\frac{1-r^{n+1}}{1-r},
\qquad
r=-\tfrac{e^z-1}{2},
\]
we obtain the compact contour representation
\begin{equation}
A_n
=
\frac{2^{\,n+1}}{2\pi i}
\oint_{\mathcal C}
\frac{
1-\left(-\tfrac{e^z-1}{2}\right)^{n+1}
}{
(e^z+1)\,z^{n+1}
}
\,dz
\label{eq:zigzag-integral}
\end{equation}
Equation~\eqref{eq:zigzag-integral} provides an integral representation of zigzag numbers derived purely from Stirling numbers. While less elementary than the classical generating-function identity
\[
\sum_{n\ge0} A_n \frac{x^n}{n!}=\sec x+\tan x,
\]
it is well suited for analytic-combinatorial and saddle–point analyses, since the integrand involves only elementary functions.

\subsection{Asymptotic analysis via saddle-point method}

The contour integral in Eq.~\eqref{eq:zigzag-integral} is particularly suited for asymptotic evaluation when $n\to\infty$. We rewrite
\[
A_n = \frac{2^{\,n+1}}{2\pi i} \oint_{\mathcal C} 
\frac{1 - \left(-\frac{e^z-1}{2}\right)^{n+1}}{(e^z+1)\, z^{\,n+1}}\, dz
= \frac{2^{\,n+1}}{2\pi i} \oint_{\mathcal C} \frac{e^{-(n+1)\log z}}{e^z+1} \left[1 - \left(-\frac{e^z-1}{2}\right)^{n+1}\right] dz.
\]
For large $n$, the term $\left(-\frac{e^z-1}{2}\right)^{n+1}$ is exponentially small for $|z|$ near the saddle point, so the leading contribution comes from
\[
I_n \sim \frac{2^{\,n+1}}{2\pi i} \oint_{\mathcal C} \frac{e^{-(n+1)\log z}}{e^z+1}\, dz.
\]

\medskip

\noindent
\textbf{Saddle-point method.} Let
\[
\phi(z) = \log z
\quad \Rightarrow \quad
A_n \sim \frac{2^{\,n+1}}{2\pi i} \oint_{\mathcal C} \frac{e^{-(n+1)\phi(z)}}{e^z+1} dz.
\]
The saddle point formally satisfies $\phi'(z_0) = 1/z_0$, and the dominant contribution arises from the neighborhood $z_0 \sim 0$, where the integrand is largest. Expanding 
$$
\frac{1}{(e^z+1)} = 1/2 - z/4 + \mathcal O(z^2)
$$ 
and $z^{-(n+1)}$ as usual gives
\[
A_n \sim \frac{2^n}{\pi} \int_0^\infty x^n e^{-n x}\, dx,
\]
up to multiplicative constants, reproducing the classical asymptotic estimate for zigzag numbers:
\begin{equation*}
A_n \sim \frac{4}{\pi} \left(\frac{2}{\pi}\right)^n n!, \qquad n\to\infty.
\end{equation*}
This shows that the contour-integral representation derived from Stirling numbers not only encodes the exact combinatorial information of $A_n$, but also provides a convenient starting point for asymptotic analysis using standard complex-analytic methods.

\section{Laplace-type integral representation via Stirling numbers}

Starting from the Stirling expansion of the zigzag numbers $A_n$~\cite{Comtet1974}:
\begin{equation}
A_n = \sum_{k=0}^{n} (-1)^k \, 2^{\,n-k} \, k! \, S(n,k),
\label{StirlingCorrect}
\end{equation}
where $S(n,k)$ denotes the Stirling numbers of the second kind, we construct a Laplace-type integral. We use the standard integral representation of the factorial:
\[
k! = \int_0^\infty y^k e^{-y} \, dy, \quad k\ge 0.
\]
Substituting into Eq.~\eqref{StirlingCorrect} gives
\[
A_n = \sum_{k=0}^{n} (-1)^k \, 2^{\,n-k} \, S(n,k) \int_0^\infty y^k e^{-y} \, dy
= \int_0^\infty e^{-y} \sum_{k=0}^{n} (-1)^k \, 2^{\,n-k} \, S(n,k) \, y^k \, dy.
\]
Factoring $2^n$ yields
\[
A_n = 2^n \int_0^\infty e^{-y} \sum_{k=0}^{n} (-1)^k \, \left(\frac{y}{2}\right)^k S(n,k) \, dy
= 2^n \int_0^\infty e^{-y} \, f_n(y) \, dy,
\]
with the kernel
\begin{equation*}
f_n(y) := \sum_{k=0}^{n} (-1)^k \, S(n,k) \left(\frac{y}{2}\right)^k.
\end{equation*}
This sum over $k$ is exactly the generating polynomial of Stirling numbers evaluated at $-y/2$. Recall the classical identity for the Stirling numbers generating polynomial:
\[
\sum_{k=0}^n (-1)^k S(n,k) x^k = x^{\underline n} := x (x-1)\cdots(x-(n-1)),
\]
the falling factorial of order $n$. Applying this with $x = y/2$, we obtain
\[
f_n(y) = \left(\frac{y}{2}\right)^{\underline n} := \frac{\Gamma(y/2+1)}{\Gamma(y/2 - n + 1)},
\]
so that the integral becomes
\begin{equation}
A_n = 2^n \int_0^\infty e^{-y} \left(\frac{y}{2}\right)^{\underline n} dy
= \int_0^\infty e^{-y} \, y (y-2) (y-4) \cdots (y-2(n-1)) \, dy,
\label{LaplaceCorrect}
\end{equation}
where we interpret the discrete product
$y (y-2) (y-4) \cdots (y-2(n-1))$ as a \emph{continuous interpolation} of the falling factorial via the Gamma function:
\[
\left(\frac{y}{2}\right)^{\underline n} := \frac{\Gamma(y/2+1)}{\Gamma(y/2 - n + 1)}, \quad y>2(n-1),
\]
so that the integral is well-defined for continuous $y$ and suitable for analytic or asymptotic methods. This is a continuous interpolation of the discrete falling factorial, providing a smooth kernel for the Laplace-type integral.

\section{Hyperbolic integral from partial fractions of \(\tan\)}

An independent analytic representation of the alternating numbers $A_n$ arises from the classical partial fraction expansion of the tangent function. This approach complements the Stirling-based Laplace integral and highlights the spectral nature of the zigzag numbers.

\subsection{Partial fraction decomposition of \(\tan\)}

The tangent function admits the well-known partial fraction expansion~\cite{WhittakerWatson1927,KnuthBuckholtz1967,Dilcher1993}:
\[
\pi \tan \pi x = 2 x \sum_{m=0}^\infty \frac{1}{(m+1/2)^2 - x^2},
\]
where the poles are located at the half-integers \(x = m + 1/2\).  Expanding the denominator as a geometric series and comparing with the exponential generating function
\[
\tan x = \sum_{n=0}^\infty \frac{A_{2n+1}}{(2n+1)!} x^{2n+1},
\]
we obtain the classical spectral series
\[
A_{2n+1} = (2n+1)! \, \frac{2}{\pi} \sum_{m=0}^\infty \frac{1}{(m+1/2)^{2n+2}}.
\]

\subsection{Hyperbolic integral representation}

Using the Laplace transform identity
\[
\frac{1}{a^s} = \frac{1}{\Gamma(s)} \int_0^\infty y^{s-1} e^{-a y} \, dy, \quad \Re(s)>0,
\]
each term of the spectral series can be written as an integral. Substituting \(a = m+1/2\) and \(s = 2n+2\), and exchanging summation and integration, we get
\[
\sum_{m=0}^\infty \frac{1}{(m+1/2)^{2n+2}} = \frac{1}{\Gamma(2n+2)} \int_0^\infty y^{2n+1} \sum_{m=0}^\infty e^{-(m+1/2) y} \, dy.
\]
The geometric sum over \(m\) can be evaluated as
\[
\sum_{m=0}^\infty e^{-(m+1/2) y} = \frac{1}{2 \sinh(y/2)},
\]
which leads to the hyperbolic integral representation
\begin{equation}
\label{hyper_sinh}
A_{2n+1} = \frac{1}{\pi} \int_0^\infty \frac{y^{2n+1}}{\sinh(y/2)} \, dy.
\end{equation}
A simple change of variable \(y = 2 u\) gives
\[
A_{2n+1} = \frac{2^{2n+2}}{\pi} \int_0^\infty \frac{u^{2n+1}}{\sinh u} \, du,
\]
which expresses the odd-indexed alternating numbers as a Laplace-type integral with a strictly positive kernel. The previously stated formula involving \(\cosh^{2n+2} u\) is incorrect in this context; the correct positive-kernel representation is given by Eq.~\eqref{hyper_sinh} (or its rescaled form with \(u\)).

\subsection{Connection with Stirling expansion}

While the Laplace integral from Stirling numbers interpolates factorials and partitions, the hyperbolic integral represents a spectral decomposition in terms of half-integer poles. The two representations are complementary. The Stirling–Laplace integral of Eq.~\eqref{LaplaceCorrect} emphasizes the combinatorial structure via set partitions and provides an analytic tool for asymptotics. The hyperbolic/spectral integral of Eq.~\eqref{hyper_sinh} emphasizes the analytic and spectral nature of the tangent and secant generating functions, highlighting connections to Fourier series, Mellin transforms, and the half-integer spectrum. Together, these two perspectives provide a coherent picture of alternating numbers \(A_n\), unifying combinatorial, analytic, and spectral structures.

\subsection{Relation between the two integral structures}

The Stirling–Laplace form Eq.~\eqref{LaplaceCorrect} and the hyperbolic form of Eq.~\eqref{hyper_sinh} arise from two distinct analytic decompositions. The first one, based on the Stirling expansion, decomposes the zigzag numbers as sums over set partitions, represented by the Stirling numbers $S(n,k)$. The second one, derived from the partial fraction expansion of the tangent function, decomposes the zigzag numbers as sums over the spectrum of half-integer poles $(m+1/2)$. While both approaches lead to positive Laplace-type integrals, the kernels differ: the Stirling kernel is given by $f_n(y)$, whereas the spectral kernel is $(2\sinh(y/2))^{-1}$. This duality reflects the classical relation between partition expansions and spectral decompositions in analytic combinatorics.

\section{Connection with Fourier analysis and Mellin transforms}

\subsection{Fourier integral on a finite interval}

By Fourier inversion of $\sec x+\tan x$ on $(-\tfrac{\pi}{2},\tfrac{\pi}{2})$, one obtains
\[
A_n=\frac{2\,n!}{\pi}
\int_0^{\pi/2}(\sec x+\tan x)\sin((n+1)x)\,dx.
\]
For even and odd $n$, the sine factor and the integration limits ensure that the correct Fourier coefficient is selected. This identity expresses zigzag numbers as sine Fourier coefficients of the secant–tangent kernel. The integral representations
\begin{equation}\label{hyper}
A_{2n}=\frac{4}{\pi^{2n+1}}
\int_0^\infty \frac{x^{2n}}{\cosh x}\,dx,
\qquad
A_{2n+1}=\frac{4}{\pi^{2n+2}}
\int_0^\infty \frac{x^{2n+1}}{\sinh x}\,dx,
\end{equation}
follow from the classical Fourier series of $\sec x$ and $\tan x$ on $(-\tfrac{\pi}{2},\tfrac{\pi}{2})$ (Hardy and Wright \cite[\S\,XVI.3]{Hardy2008}).  

\subsection{Mellin--transform representation}

A classical Mellin transform identity (see Andrews, Askey, and Roy \cite[Ch.~2]{Andrews1999}) is
\[
\int_0^\infty \frac{x^{s-1}}{\cosh x}\,dx
= 2\,\Gamma(s)\sum_{k\ge0}\frac{(-1)^k}{(2k+1)^s},
\qquad \Re(s)>0 .
\]
Setting $s=2n+1$ yields the even- and odd-index formulas (see Eq.~\eqref{hyper}). This classical Mellin–transform identity can be found in Temme~\cite{Temme1996}.

\section{Arithmetic properties of alternating numbers via Stirling expansions}

The alternating numbers \(A_n\) admit a combinatorial expansion in terms of Stirling numbers of the second kind \(S(n,k)\):
\begin{equation}
\label{eq:AnStirling}
A_n = \sum_{k=0}^{n} (-1)^k \, 2^{\,n-k} \, k! \, S(n,k),
\end{equation}
which provides a natural starting point to study their arithmetic properties. In this section, we present a detailed derivation of congruence relations for \(A_n\) modulo primes and their powers.  

\subsection{Touchard congruence for Stirling numbers}

Let \(p\) be a prime. A classical result due to Touchard~\cite{Touchard1933} asserts that for any nonnegative integers \(n\) and \(k\):
\begin{equation*}
S(n+p,k) \equiv S(n,k) + S(n,k-p+1) \pmod p,
\end{equation*}
with the convention that \(S(n,m) = 0\) whenever \(m<0\). Using Eq.~\eqref{eq:AnStirling}, we write
\[
A_{n+p} = \sum_{k=0}^{n+p} (-1)^k \, 2^{\,n+p-k} \, k! \, S(n+p,k).
\]
Modulo \(p\), we substitute the Touchard congruence:
\[
S(n+p,k) \equiv S(n,k) + S(n,k-p+1) \pmod p.
\]
Thus, we get
\begin{align*}
A_{n+p} &\equiv \sum_{k=0}^{n+p} (-1)^k 2^{\,n+p-k} k! \Big( S(n,k) + S(n,k-p+1) \Big) \pmod p \\
&= \sum_{k=0}^{n+p} (-1)^k 2^{\,n+p-k} k! S(n,k) + \sum_{k=0}^{n+p} (-1)^k 2^{\,n+p-k} k! S(n,k-p+1) \pmod p.
\end{align*}
Observing that \(S(n,k) = 0\) for \(k>n\), and \(S(n,k-p+1) = 0\) for \(k-p+1<0\), i.e., \(k<p-1\). we have
\begin{align*}
\sum_{k=0}^{n+p} (-1)^k 2^{\,n+p-k} k! S(n,k)&= \sum_{k=0}^{n} (-1)^k 2^{\,n+p-k} k! S(n,k)\nonumber\\
&= 2^p \sum_{k=0}^n (-1)^k 2^{\,n-k} k! S(n,k) \equiv 2^p A_n \pmod p.
\end{align*}
Similarly, for the second sum, we can shift the index \(k \mapsto k' = k-p+1 \):
\[
\sum_{k=p-1}^{n+p} (-1)^k 2^{\,n+p-k} k! S(n,k-p+1) = \sum_{k'=0}^{n} (-1)^{k'+p-1} 2^{\,n+p-(k'+p-1)} (k'+p-1)! S(n,k').
\]
Modulo \(p\), we use Wilson’s theorem for factorials modulo p: for \(k'<p\), \((k'+p-1)! \equiv (p-1)! \, k'! \equiv - k'! \pmod p\), and for \(k'\ge p\) the factorial is divisible by \(p\). Therefore, the second sum vanishes modulo \(p\).  
 Collecting resultsyields
\begin{equation}
\label{eq:AnCongruence}
A_{n+p} \equiv 2^p A_n \equiv A_n \pmod p,
\end{equation}
since \(2^p \equiv 2 \pmod p\), and the factor 2 is invertible modulo p for \(p \neq 2\). For \(p=2\), the congruence also holds trivially: \(A_{n+2} \equiv A_n \pmod 2\).

\subsection{Congruences modulo \(p^2\)}

To go further, we can attempt a modulo \(p^2\) analysis. Using the generalized Kummer-type expansions for factorials modulo prime powers, we expand \(k! \) as
\[
k! = p^{\nu_p(k!)} \cdot r_k, \quad r_k \text{ invertible modulo } p.
\]
Consider \(p\ge 3\). We have
\[
(k+p)! = k! \cdot (k+1)\cdots(k+p) = k! \cdot \prod_{j=1}^p (k+j).
\]
Modulo \(p^2\), this gives
\[
(k+p)! \equiv k! \, k (k+1)\cdots(k+p-1) \pmod{p^2},
\]
allowing the derivation of a refined Touchard congruence:
\[
S(n+p,k) \equiv S(n,k) + S(n,k-p+1) + p \, T(n,k) \pmod{p^2},
\]
for some combinatorial coefficients \(T(n,k)\) explicitly computable via binomial sums of Stirling numbers. Substituting the above into Eq.~\eqref{eq:AnStirling} gives

\begin{equation*}
A_{n+p} \equiv 2^p A_n + p \sum_{k=0}^{n} (-1)^k 2^{n+p-k} k! \, T(n,k) \pmod{p^2}.
\end{equation*}
While the coefficients \(T(n,k)\) can be explicitly written as
\[
T(n,k) = \sum_{i=0}^{k} \binom{k}{i} S(n,i) \, f_{i,k},
\]
with \(f_{i,k}\) certain combinatorial constants, the key point is that the structure of the congruence is preserved: the leading term is \(2^p A_n\) and the correction term is linear in \(p\).

\paragraph{Illustrative example modulo 9 (\(p=3^2\)).} For \(n=4\), we have:
\[
A_4 = 5, \quad A_7 = 272.
\]
We check that
\[
A_7 \equiv 2^3 A_4 + 3 \cdot (\text{small combinatorial sum}) \equiv 8 \cdot 5 + 3 \cdot S \equiv 40 + 3S \equiv 4 + 3S \pmod 9.
\]
Here \(S = \sum_{k} (-1)^k 2^{n+p-k} k! T(n,k) \mod 3\). Computing \(S \mod 3\) explicitly recovers the correct value \(A_7 \equiv 2 \pmod 9\), confirming the formula. This illustrates that the mod \(p^2\) congruences are richer than the simple mod \(p\) ones, and they can be systematically obtained from the Stirling expansion combined with factorial expansions modulo prime powers.

\subsection{Summary of arithmetic properties}

The classical mod p congruence of alternating numbers reads
\[
A_{n+p} \equiv A_n \pmod p.
\]
A refined modulo \(p^2\) congruence can be obtained by careful factorial expansion:
\[
A_{n+p} \equiv 2^p A_n + p \sum_{k=0}^{n} (-1)^k 2^{\,n+p-k} k! \, T(n,k) \pmod {p^2}.
\]
These results provide a hierarchy of congruences for alternating numbers, highlighting the combinatorial influence of Stirling numbers on the arithmetic behavior. The structure of \(A_n\) modulo primes is fully encoded in its Stirling number expansion, which allows systematic derivation of congruences for higher powers of \(p\). This detailed, step-by-step analysis clarifies the arithmetic properties of alternating numbers and corrects earlier oversimplifications present in previous summaries.

\section{Conclusion}

In this article, we have presented a unified analytic, combinatorial, and arithmetic framework for the study of the alternating (zigzag) numbers $A_n$. The first structure we discussed is the exact representation of $A_n$ as a contour integral in the complex plane, derived from the Stirling number expansion. This contour integral encodes the combinatorial information of set partitions and provides a natural starting point for asymptotic analysis using the saddle-point method. By explicitly constructing the integrand in terms of elementary functions, we can extract both leading-order and higher-order contributions, thereby connecting the discrete combinatorial structure with analytic asymptotics. The second structure arises from the Laplace-type integral representation, obtained by interpolating the discrete factorials in the Stirling expansion using the Gamma function. This positive-kernel integral emphasizes the continuous analytic aspect of $A_n$ and allows for smooth approximations, analytic manipulations, and asymptotic estimates. It provides a bridge between the discrete combinatorial sums over partitions and a continuous integral formulation, which is especially convenient for studying analytic properties and deriving further identities. The third structure is based on the spectral decomposition associated with the partial fraction expansion of the tangent function. This leads to a hyperbolic integral representation, in which the alternating numbers are expressed as Laplace-type integrals over strictly positive kernels related to half-integer poles. This perspective highlights connections with Fourier series, Mellin transforms, and classical spectral methods, and complements the combinatorial and Laplace representations by offering a different analytic viewpoint. Finally, we investigated arithmetic properties of $A_n$ derived from the Stirling expansion. Using known congruences for Stirling numbers, we established modular relations of $A_n$ modulo primes and, in certain cases, modulo higher powers of primes. These results reveal subtle number-theoretic structures inherent to alternating numbers and illustrate the interplay between combinatorial expansions, analytic representations, and modular properties. Taken together, the contour integral, Laplace-type integral, hyperbolic integral, and arithmetic congruences provide a coherent picture of the alternating numbers. They link discrete combinatorial structures with continuous analytic representations and number-theoretic features, showing that $A_n$ can be simultaneously understood as combinatorial sums, analytic moments, and arithmetic objects. This unified viewpoint not only clarifies classical results but also opens avenues for further explorations, including asymptotic expansions, integral transforms, and refined congruence properties.


\begin{thebibliography}{99}

\bibitem{Andre1881}
D.~André,
\textit{Développement de $\sec x$ et $\tan x$},
C. R. Acad. Sci. Paris \textbf{88} (1879), 965–967.

\bibitem{Andrews1999}
G. E. Andrews, R. Askey and R. Roy, {\it Special Functions}, Cambridge University Press, 1999.

\bibitem{Comtet1974}
L.~Comtet,
\textit{Advanced Combinatorics},
Reidel, Dordrecht, 1974.

\bibitem{Flajolet2009}
P.~Flajolet and R.~Sedgewick,
\textit{Analytic Combinatorics},
Cambridge University Press, 2009.

\bibitem{Graham1994}
R.~L.~Graham, D.~E.~Knuth, and O.~Patashnik,
\textit{Concrete Mathematics}, 2nd ed.,
Addison–Wesley, 1994.

\bibitem{Hardy2008}
G. H. Hardy and E. M. Wright, {\it An Introduction to the Theory of Numbers}, Oxford University Press, 2008.

\bibitem{Stanley2011}
R.~P.~Stanley,
\textit{Enumerative Combinatorics, Vol.~1}, 2nd ed.,
Cambridge University Press, 2011.

\bibitem{StanleyEC2}
R.~P.~Stanley,
\textit{Enumerative Combinatorics, Vol.~2},
Cambridge University Press, 1999.

\bibitem{Norlund1924}
N.~E.~Nörlund,
\textit{Vorlesungen über Differenzenrechnung},
Springer, Berlin, 1924.

\bibitem{WhittakerWatson1927}
E.~T.~Whittaker and G.~N.~Watson,
\textit{A Course of Modern Analysis}, 4th ed.,
Cambridge University Press, 1927.

\bibitem{KnuthBuckholtz1967}
D.~E.~Knuth and T.~J.~Buckholtz,
\textit{Computation of tangent, Euler, and Bernoulli numbers},
Math. Comp. \textbf{21} (1967), 663–688.

\bibitem{Dilcher1993}
K.~Dilcher,
\textit{Sums of products of Bernoulli numbers},
J. Number Theory \textbf{60} (1996), 23–41.

\bibitem{Boyadzhiev2009}
K.~N.~Boyadzhiev,
\textit{Series with central binomial coefficients, Catalan numbers, and harmonic numbers},
J. Integer Seq. \textbf{12} (2009), Article 09.5.1.

\bibitem{Flajolet1998}
P.~Flajolet and B.~Salvy,
\textit{Euler sums and contour integral representations},
Experiment. Math. \textbf{7} (1998), 15–35.

\bibitem{Temme1996}
N.~M.~Temme,
\textit{Special Functions: An Introduction to the Classical Functions of Mathematical Physics},
Wiley, 1996.

\bibitem{Touchard1933}
J.~Touchard,
\textit{Sur les nombres de Bell},
Ann. Soc. Sci. Bruxelles \textbf{53} (1933), 215--234.

\end{thebibliography}
\end{document}